\documentclass[10pt,twoside]{siamltex}
\usepackage{amsfonts,epsfig}

\begin{document}

\newcommand{\lanbox}{\hfill \hbox{$\, 
\vrule height 0.25cm width 0.25cm depth 0.01cm
\,$}}

\bibliographystyle{plain}
\title{Short proofs of theorems of  Mirsky and Horn on diagonals and eigenvalues of matrices}
\author{ Eric A. Carlen\thanks{Department of Mathematics, Hill Center,
Rutgers University,
110 Frelinghuysen Road,
Piscataway, NJ 08854-8019. Work partially
supported by U.S. National Science Foundation
grant DMS 06-00037. }
\and
Elliott H. Lieb\thanks
 Departments of Mathematics and Physics, Jadwin
Hall, Princeton University, P.~O.~Box 708, Princeton, NJ
  08542. Work partially
supported by U.S. National Science Foundation
grant PHY 06-52854.
\hfill\break
\copyright\, 2009 by the authors. This paper may be reproduced, in its
entirety, for non-commercial purposes.
}
\maketitle
\thispagestyle{empty}
\pagestyle{empty}

\begin{abstract}
A theorem of Mirsky provides  necessary and sufficient conditions for the existence of an $N$-square
complex matrix with prescribed diagonal entries and prescribed eigenvalues. We give a simple inductive proof of this theorem. 
\if false
.\\

\medskip

\centerline{\bf July 3, 2009}
\fi
\end{abstract}

\begin{keywords}
Key Words: majorization, eigenvalues,
prescribed diagonals
\end{keywords}
\begin{AMS}
15A42, 15A51.
\end{AMS}

\renewcommand{\thetheorem}{1.\arabic{equation}}

\renewcommand{\theequation}{1.\arabic{equation}}

Theorems of Alfred Horn \cite{horn} and Leon Mirsky \cite{mirsky} 
give necessary and sufficient conditions for the existence of an $N$-square matrix $A$  that has prescribed  diagonal elements and prescribed eigenvalues. In the case of Horn's Theorem,
$A$ is required to be Hermitian, and in the case of Mirsky's theorem no constraint is  imposed on $A$.
While  simple inductive 
proofs of Horn's Theorem  have been published, the only published
proof of Mirsky's Theorem that we are aware of is the original, which is neither simple nor transparent.
In this note, we give a simple proof of Mirsky's Theorem
which we relate to 
 Chan and Li's  proof  \cite{chanli} of Horn's Theorem. We begin by recalling the results in \cite{horn,mirsky} in more detail.

Let $\Lambda =(\lambda_1, \dots, \lambda_N)$ and 
$D=(d_1,\dots, d_N) $ be two sequences of $N$ real numbers.  Recall that  $\Lambda$  majorizes  $D$  (we write  $\Lambda \succ D$) if
for each $k =1,\dots,N$, the sum of the $k$ largest $\lambda$'s is at least as
large as the sum of the $k$ largest $d$'s, with equality for $k=N$. 
If $\Lambda$ and $D$ are the eigenvalue and diagonal sequences,
respectively, of an $N$-square  Hermitian matrix $A$, then $A=U[ \Lambda] U^*$ for some
unitary $U$, where  $[\Lambda]$  denotes the diagonal matrix with
diagonal sequence $\Lambda$.  Thus,
\begin{equation}\label{rel}
{\textstyle d_i = \sum_{j=1}^N |U_{ij}|^2 \lambda_j\qquad{ \rm for}\quad 1 \leq i \leq N\ .}
\end{equation}

These $N$ equations
lead directly to
the necessary  condition $\Lambda \succ D$  found by Schur in 1923 \cite{schur}.
Horn's Theorem \cite{horn}  says that this necessary condition is also sufficient,
and moreover, $U$ can be taken to be real orthogonal.
Thus, while Birkhoff showed \cite{birkhoff} that $\Lambda \succ D $ if and only if
there exists a bi-stochastic matrix $S$ such that $D=S\Lambda $, Horn's
Theorem states that $S$ can always be chosen to lie in the
smaller set of {\it ortho-stochastic} matrices.

Mirsky \cite{mirsky} gave another proof of Horn's Theorem, and  also  proved the following:
Given two complex $N$-sequences $\Lambda$ and $D$,  there exists
an $N$-square matrix whose eigenvalue sequence is $\Lambda$ and whose diagonal sequence is $D$ 
if and only if
\begin{equation}\label{cond1}
{\textstyle \sum_{j=1}^N\lambda_j = \sum_{j=1}^N d_j\ ,}
\end{equation}
and that if $\Lambda$ and $D$ are both real, then the matrix can be taken to be real as well.

We now give a simple inductive 
proof of Mirsky's Theorem. 
Let $[T_\Lambda]$ denote the $N$-square matrix that has
$\Lambda$ as its diagonal sequence, and has $1$ in every entry immediately above the diagonal, and
$0$ in all remaining entries. A {\em unit lower triangular} matrix, is a square matrix in which every entry on the diagonal is $1$, and every entry above the diagonal is $0$. 
The unit lower triangular matrices are a group under matrix multiplication. 
The following result includes Mirsky's Theorem and a little more.

\smallskip

\begin{theorem}\label{mirtheorem}
There exists a matrix $A$ with eigenvalue sequence $\Lambda$ and diagonal
sequence $D$ if and only if condition (\ref{cond1}) is satisfied.  Indeed,
under this condition, there exists a unit lower triangular
matrix $L$ such that the diagonal sequence of $L^{-1}[T_\Lambda]L$ is $D$.
If $\Lambda$ and $D$ are both real, $L$ can be taken to be real.
\end{theorem}
\smallskip

\noindent{\bf Proof:} The necessity  of   (\ref{cond1}) for the existence of 
$A$ is obvious, and we must prove that (\ref{cond1}) is sufficient for the existence of $L$.
For $N=2$,
direct computation shows that, with $L =  \left[\begin{array}{cc}
1 &0  \\
c & 1  
 \end{array}
\right]$,
the diagonal sequence of $L^{-1}[T_\Lambda]L$ is $(\lambda_1+c,\lambda_2-c)$.
Choosing $c= d_1-\lambda_1$ yields the assertion for $N=2$.

Assume that the assertion is proved for $K$-square matrices of size up through $K = N-1$.  Given  $N$-sequences $\Lambda$ and $D$
satisfying (\ref{cond1}), let $L$ be the $2$-square unit lower triangular matrix such that
 $L^{-1}[T_{(\lambda_1,\lambda_2)}]L$ has the diagonal sequence $(d_1,\widetilde \lambda_2)$, where
 \begin{equation}\label{cond8}
\widetilde  \lambda_2 = \lambda_1+\lambda_2 - d_1\ .
 \end{equation}
 Let $\widetilde L$ be the $N$-square
 matrix obtained by replacing the upper left $2$-square block of $I_{N\times N}$ by $L$.
 Then $\widetilde L^{-1}[T_\Lambda]\widetilde L$ has the following properties: Its upper left
 $2$-square block has diagonal sequence $(d_1,\widetilde \lambda_2)$, its lower right $(N-2)$-square
 block is the same as that  of $[T_\Lambda]$, and 
$\left[\widetilde L^{-1}[T_\Lambda]\widetilde L\right]_{2,3} = 1$. 
 
 The lower right $(N-1)$-square block of   $\widetilde L^{-1}[T_\Lambda]\widetilde L$
 is $[T_{\widetilde \Lambda}]$ where $\widetilde \Lambda$ is the $(N-1)$-sequence
 $(\widetilde  \lambda_2, \lambda_3,\dots,\lambda_N)$.  By the inductive hypothesis, there exists
  a unit lower triangular matrix $M$ such that $M^{-1}[T_{\widetilde \Lambda}]M$ has the diagonal sequence $(d_2,\dots,d_N)$. Now let $\widetilde M$ be the $N$-square matrix obtained by
  replacing the lower right $(N-1)$ square block of $I_{N\times N}$ with $M$. Then
 $(\widetilde L\widetilde M)^{-1}[T_\Lambda](\widetilde L\widetilde M)$
  has the diagonal sequence $D$, and  $\widetilde L\widetilde M$ is unit lower triangular.
\lanbox

Our proof of Theorem~\ref{mirtheorem} is related to a proof of Horn's Theorem  by Chan and Li \cite{chanli}. For the reader's convenience, we briefly recall their proof  in our notation, and then conclude with a brief comparison of the proofs.

\smallskip

\noindent{\bf Proof of Horn's Theorem, after Chan and Li:} 
Assume that $\Lambda$ and $D$ are in decreasing order. 
If $N=2$, then $\Lambda \succ D $ implies 
$\lambda_1 \ge d_1 \ge d_2 \ge \lambda_2$. For $\lambda_1>\lambda_2$,
$$U =  (\lambda_1 - \lambda_2)^{-1/2}
 \left[\begin{array}{cc}
(d_{1}-\lambda_{2}) ^{1/2} &-(\lambda_{1}-d_{1})^{1/2}  \\
(\lambda_{1}-d_{1})^{1/2} & \phantom{-}(d_{1}-\lambda_{2}) ^{1/2} 
 \end{array}
\right]$$
is a real orthogonal matrix, and the diagonal sequence of $U^*[\Lambda]U$ is $D = (d_1,d_2)$.
For $\lambda_1 = \lambda_2$, $\Lambda =D$,  and we may take $U = I$.  This proves Horn's Theorem for $N=2$. 

We proceed inductively: Let $\Lambda$ and $D$ be two real $N$-sequences with 
$\Lambda \succ D $.
There is a $K$ with  $1\leq K < N$ such that $\lambda_{K} \geq d_{K} \geq
d_{K+1} \geq \lambda_{K+1}$: Take $K$ to be the smallest $j<N$
such that $d_{j+1} \geq \lambda_{j+1}$. 
Note that $(\lambda_{K}, \lambda_{K+1}) \succ (d_{K},\lambda_{K+1}')$, where
\begin{equation}\label{rel2}
\lambda_{K+1}' = \lambda_{K+1}+\lambda_{K} -d_{K} \ .
\end{equation}
Our argument in the case $N=2$ shows that there exists a $2$-square real orthogonal matrix $U$
such that $U^*[(\lambda_K,\lambda_{K+1})]U$ has $(d_{K},\lambda_{K+1}')$
as its diagonal sequence. 

We use $U$ to construct an $N$-square real orthogonal matrix $T$
as follows:   Start with the identity matrix $I_{N\times N}$, and replace the $2$-square
diagonal block at positions $K,K+1$ by the $2$-square matrix $U$.   Then
 $T^*[\Lambda]T$ has the diagonal sequence 
$\Lambda' $, which is obtained from $\Lambda$ 
by replacing $\lambda_{K}$ by $d_{K}$ and
$\lambda_{K+1}$ by $\lambda_{K+1}'$ in $\Lambda$. We note the important fact that
$\Lambda' \succ D$.

What follows is especially
simple if $K=1$, and the reader may wish to consider that case first. 
Let $\Lambda''$ and $D''$ be the
sequences of $N-1$ real numbers obtained by deleting $d_{K}$, the common $K$th term in
both $\Lambda'$ and $D$. Then $\Lambda''\succ D''$, and hence, by
induction, there exists an $(N-1)$-square orthogonal matrix $V$ such
that $V[\Lambda'' ]V^*$ has the diagonal sequence $D''$.  It is
convenient to index the entries of $V$ using $\{\dots,K-1,K+1,\dots\}$
leaving out the ``deleted'' index $K$. Promote $V$ to an
$N$-square orthogonal matrix by setting $V_{K,K} =1$ and $V_{j,K} =
V_{K,j} = 0$ for $j\neq K$. Then $(VT)[\Lambda](VT)^*$ has $D$ as
its diagonal sequence.  
 \lanbox
\medskip

While Chan and Li's proof of Horn's Theorem and our  proof of Mirsky's Theorem 
have a similar structure, there is a significant difference: 
The matrices whose existence is guaranteed by Horn's Theorem are all necessarily
similar to the diagonal matrix $[\Lambda]$, so that $[\Lambda]$ provides a suitable starting point for
the proof.  This is not the case for Mirsky's Theorem: Indeed, if
all  entries of $\Lambda$ are equal, $[\Lambda]$ is a multiple of the identity. Thus, the set of all matrices
similar to $[\Lambda]$ contains only $[\Lambda]$  itself. Except in the  case $D = \Lambda$,
it does not contain any matrices with diagonal sequence $D$.  However, our proof shows that in the set of all
matrices similar to $[T_\Lambda]$, there is always one with diagonal sequence $D$ provided that  (\ref{cond1})
is satisfied.

\medskip

\noindent{\bf Acknowledgement:} 
We thank Prof. Roger Horn for invaluable advice.

\end{document}